\documentclass[12pt]{amsart}
\usepackage{amscd, amsfonts, amssymb}
\usepackage{verbatim}

\newcommand\CA{{\mathcal A}}
\newcommand\CB{{\mathcal B}}

\newcommand\CI{{\mathcal I}} 
\newcommand\CN{{\mathcal N}} 

\newcommand\CP{{\mathcal P}}

\newcommand\CR{{\mathcal R}}

\newcommand\fb{{\mathfrak b}}
\newcommand\fg{{\mathfrak g}}
\newcommand\fh{{\mathfrak h}}

\newcommand\fu{{\mathfrak u}}

\newcommand\fn{{\mathfrak n}}

\newcommand\cor{{\operatorname{cr}}} 

\newcommand\Lie{{\operatorname{Lie\,}}}

\newcommand\Ad{{\operatorname{Ad\,}}}

\newcommand\nil{{\operatorname{nil}}} 

\numberwithin{equation}{section}

\newtheorem{proposition}[equation]{Proposition}

\newtheorem{remark}[equation]{Remark}

\title[Alperin's Conjecture for Algebraic Groups]
{Alperin's Conjecture for Algebraic Groups}

\author[G. R\"ohrle]{Gerhard R\"ohrle}
\address
{Fakult\"at f\"ur Mathematik,
Ruhr-Universit\"at Bochum,
Universit\"ats\-strasse 150,
D-44780 Bochum, Germany}
\email{gerhard.roehrle@rub.de}

\author[R. Rouquier]{Rapha\"el Rouquier}
\address
{Mathematics Institute,
University of Oxford, 
24-29 St Giles', Oxford, OX1 3LB, United Kingdom}
\email{rouquier@maths.ox.ac.uk}


\makeatother

\newcommand{\Size}[1]{\left|#1\right|}

\begin{document}

\begin{abstract}
We prove analogues for reductive algebraic groups of some results 
for finite groups due to Kn\"orr and Robinson from \cite{knoerrrobinson} 
which play a central r\^ole in their reformulation of Alperin's conjecture
for finite groups.
\end{abstract}

\maketitle

\section{Introduction}
\label{sec:intro}

Let $G$ be a finite group, $p$ a prime and $k$ an algebraically 
closed field of characteristic $p$.
By $kG$ we denote the modular group algebra of $G$.
Alperin's conjecture \cite{alperin} 
asserts that the number of isomorphism classes of simple $kG$-modules 
equals the sum of the number of 
isomorphism classes of projective simple $k[N_G(P)/P]$-modules 
where $P$ is a $p$-subgroup of $G$ and the sum is taken over all 
$p$-subgroups $P$ of $G$ up to $G$-conjugacy.
Kn\"orr and Robinson \cite[Thm.\ 3.8]{knoerrrobinson}
reformulated this conjecture 
in terms of the vanishing of an alternating sum of 
the number of simple modules for normalizers of $p$-subgroups. 
More precisely,
they showed that Alperin's conjecture holds for all finite groups if and
only if their alternating sum conjecture holds for all finite groups.
For finite groups of Lie type, Alperin's original conjecture was 
first proved by M.\ Cabanes \cite{cabanes}, see also 
\cite[Thm.\ 5.3]{knoerrrobinson}, 
\cite{lehrerthevenaz}, and \cite{thevenazwebb}. 

The aim of this note it to prove analogues 
for reductive algebraic groups of some results 
of Kn\"orr and Robinson from \cite{knoerrrobinson} that are relevant in
their reformulation of Alperin's conjecture.


\section{Complexes of nilpotent subalgebras of $\fg$}
\label{sec:complex}


Let $G$ be a connected reductive linear algebraic group
defined over an algebraically closed field $k$. 
We denote the Lie algebra of $G$ by $\Lie G$ or by $\fg$; 
likewise for closed subgroups of $G$.
For a closed subgroup $H$ of $G$, 
the 
normalizer of $\Lie H = \fh$ in $G$ is defined by 
$N_G(\fh) = \{g\in G \mid \Ad g(\fh) \subseteq \fh \}$, 
where $\Ad g$ denotes the adjoint action of $g\in G$ on $\fg$.

By $R_u(H)$ we denote the unipotent radical of $H$ and frequently write
$\nil(\fh)$ for the nilradical $\Lie(R_u(H))$ of $\fh$.

We define several simplicial complexes consisting of various 
chains of nilpotent subalgebras of $\fg$. They are analogues of the 
subcomplexes of $p$-subgroups in finite group theory mentioned above.
To our knowledge they have not been studied yet in the context of 
reductive algebraic groups.

Let $\CN$ denote the simplicial complex associated to the
partially ordered set of all 
chains of nilpotent subalgebras of $\fg$.
We define $\CI$ to be the subcomplex of $\CN$ where for a fixed chain 
$C$ in $\CI$ there exists a Borel subalgebra $\fb$ of $\fg$ such that 
each member of $C$ is an ideal of $\fb$; equivalently, 
there exists a Borel subgroup $B$ of $G$ such that each member of 
$C$ is a $B$-submodule of $\nil(\Lie B)$. 
Moreover, $\CA$ is the subcomplex of $\CI$ where each member of a 
given chain $C$ is an abelian ideal of a Borel subalgebra associated 
to $C$.
Finally, by $\CR$ we denote the subcomplex of $\CI$ of chains $C$
where each member $\fn$ in $C$ satisfies 
$\fn = \nil(\Lie N_G(\fn))$.

The empty chain is considered to be a $(-1)$-simplex in each case.
We will assume that every non-empty chain $C$
in $\CN$ 
considered is of the form $\fn_0 \subset \fn_1 \subset \cdots \subset \fn_n$
where $\fn_0 = \{0\}$.
The \emph{chain stabilizer} $G_C$ 
of $C$ in $G$ is defined 
to be $G_C := \cap_{i=0}^n N_G(\fn_i)$.
We define the \emph{length} of the chain $C$ in $\CN$ by
$\Size{C}= n$, so that $\Size{C} = \dim C + 1$, where $\dim C$ is 
the dimension of $C$ as a simplex.

The adjoint representation of $G$ on $\fg$ 
induces an action of $G$ on each of the simplicial complexes defined;
for $C$ as above and $g \in G$ we define $g\cdot C$ to be the chain 
$\{0\} = (\Ad g)\fn_0 \subset (\Ad g)\fn_1 \subset \cdots \subset (\Ad g)\fn_n$. 
Let $\CN/G$ denote the set of $G$-conjugacy classes of 
chains in $\CN$; likewise for the other complexes.
Since all the chains we consider consist of nilpotent subalgebras
of $\fg$,  
we may assume that, up to $G$-conjugacy, any given  chain lies 
in the nilradical
$\nil(\fb)$ of a fixed Borel subalgebra $\fb$ of $\fg$. 
Thus, in particular, each of the sets of $G$-classes
$\CN/G$, $\CI/G$, $\CR/G$ and
$\CA/G$ is finite.

The following is the analogue of 
\cite[Prop.\ 3.3]{knoerrrobinson} 
in the context of reductive algebraic groups. 

\begin{proposition}
\label{prop:1}
Let $A$ be an abelian group and let $f$ be a 
$G$-equivariant 
function from the set of subgroups of $G$ to $A$ 
(i.e., $f$ is constant on conjugacy classes of subgroups of $G$). 
Then
\begin{align*}
& \sum_{C\in \CN/G}(-1)^{\Size{C}} f(G_C) 
 = \sum_{C\in \CI/G}(-1)^{\Size{C}} f(G_C) \\
& = \sum_{C\in \CR/G}(-1)^{\Size{C}} f(G_C) 
  = \sum_{C\in \CA/G}(-1)^{\Size{C}} f(G_C).
\end{align*}
\end{proposition}

\begin{proof}
Observe that by the remark above, each of the sums is finite.
We imitate the proof of \cite[Prop.\ 3.3]{knoerrrobinson}: 
The idea is to pair up chains which lie outside $\CA$, 
respectively outside $\CR$, so that their 
contributions in the alternating sums above cancel each other out.

First we show that the $G$-classes of 
chains in $\CI\setminus \CA$ do not contribute to 
the alternating sum 
$\sum_{C\in \CI/G}(-1)^{\Size{C}} f(G_C)$.

Let $C$ be a chain $\fn_0 \subset \fn_1 \subset \cdots \subset \fn_n$
in $\CI\setminus \CA$. Let $\fb$ be a Borel subalgebra of $G$ so that
$\fn_n \subseteq \fu = \nil(\fb)$.
We pair $C$ with a chain $C'$ in $\CI\setminus \CA$ as follows.
Since $\fn_n$ is not abelian, the commutator subalgebra
$[\fn_n,\fn_n]$ is non-trivial.
Let $j > 0$  be minimal so that 
$[\fn_n,\fn_n] \subseteq \fn_j$.
Observe that $\tilde \fn_j := \fn_{j-1} + [\fn_n,\fn_n]$
is again a $B$-submodule of $\fu$.
Now, if $\tilde \fn_j \ne \fn_j$, then we insert
$\tilde \fn_j$ (between $\fn_{j-1}$ and $\fn_j$) into $C$ 
to obtain $C'$, and 
if $\tilde \fn_j = \fn_j$, then we remove $\fn_j$ from $C$ to obtain $C'$.
In any case $C'$ again belongs to $\CI\setminus \CA$, 
since $\fn_n$ still belongs to $C'$; for, 
if $\fn_{n-1} + [\fn_n,\fn_n] = \fn_n$, then 
we have $\fn_{n-1} = \fn_n$, a contradiction. 
One readily checks that $(C')' = C$, $\Size{C'} = \Size{C} \pm 1$,
and that $(g\cdot C)' = g\cdot (C')$.
It follows that the chain stabilizers $G_C$ and $G_{C'}$ coincide. 
We may pair the contributions of the $G$-orbits of $C$ and $C'$ and
this shows that  
\[
\sum_{C\in \CI/G}(-1)^{\Size{C}} f(G_C) =
\sum_{C\in \CA/G}(-1)^{\Size{C}} f(G_C).
\]

The very same argument as the one above, 
with $C$ taken from $\CN$ instead of $\CI$, shows that in fact
the $G$-classes of 
chains in $\CN\setminus \CA$ do not contribute to 
the alternating sum 
$\sum_{C\in \CN/G}(-1)^{\Size{C}} f(G_C)$. Thus we obtain
\[
\sum_{C\in \CN/G}(-1)^{\Size{C}} f(G_C) =
\sum_{C\in \CA/G}(-1)^{\Size{C}} f(G_C).
\]

Finally, we show that chains in $\CI\setminus \CR$ do not make a contribution 
to  
$\sum_{C\in \CI/G}(-1)^{\Size{C}} f(G_C)$.
Let $C$ be a chain $\fn_0 \subset \fn_1 \subset \cdots \subset \fn_n$
in $\CI\setminus \CR$. 
We pair $C$ with a chain $C'$ in $\CI\setminus \CR$ as follows.
Let $i > 0$  be minimal so that
${\tilde \fn}_i := \nil(\Lie N_G(\fn_i)) \ne \fn_i$. 
Since $\fn_i \subseteq {\tilde \fn}_i$, we have 
${\tilde \fn}_i \not\subseteq \fn_i$, by hypothesis on $i$.
Now let $j \ge i$ be maximal so that ${\tilde \fn}_i \not\subseteq \fn_j$.
If $j < n$ and 
${\tilde \fn}_i + \fn_j = \fn_{j+1}$ we remove $\fn_{j+1}$ from $C$ and
if ${\tilde \fn}_i + \fn_j \ne \fn_{j+1}$ or if $j=n$ 
we insert ${\tilde \fn}_i + \fn_j$ into $C$ to obtain $C'$.
In any event, $C'$ still belongs to $\CI\setminus \CR$,
since $\fn_i$ is still a member of the resulting chain $C'$.  
Moreover, one readily checks that $(C')' = C$, $\Size{C'} = \Size{C} \pm 1$,
and $(g\cdot C)' = g\cdot (C')$.
Since $N_G(\fn_i) = N_G({\tilde \fn}_i)$, we obtain 
$G_C = G_{C'}$ in any case, and we 
may pair the contributions of the $G$-orbits of $C$ and $C'$;
this shows that  
\[
\sum_{C\in \CI/G}(-1)^{\Size{C}} f(G_C) =
\sum_{C\in \CR/G}(-1)^{\Size{C}} f(G_C).
\]
The result now follows.
\end{proof}

Let $\fn$ be a member of a chain $C$ in $\CR$.
Then as $\fn$ is normalized by a Borel subgroup of $G$, its 
normalizer $N_G(\fn)$ is therefore a parabolic subgroup of $G$.
Thus by definition, each member of a chain $C$ in $\CR$
is the nilradical of a parabolic subalgebra of $\fg$.
Consequently, the chain stabilizer $G_C$ is 
simply the parabolic subgroup whose nilpotent radical is the largest
member in $C$.

We define another complex of chains, $\CP$, consisting of chains
of parabolic subgroups of $G$.
If $C$ is a chain in $\CR$, then we can associate to it a chain 
$D$ in $\CP$ 
of the corresponding parabolic subgroups in $G$, that
are the normalizers  
of the members of $C$, and conversely for a chain $D$ in $\CP$ 
we can form a chain $C$ in $\CR$ by taking the nilradicals of 
the parabolic subgroups in $D$. 
Note that both operations are maps of 
complexes and both are order-reversing and preserve the lengths of chains.
Further, since parabolic subgroups are self-normalizing, 
the chain stabilizer $G_D$ of a chain $D$ in $\CP$
is simply the smallest parabolic subgroup in $D$.
In particular, if $C$ in $\CR$ and $D$ in $\CP$ correspond in this
way, then $G_C = G_D$. 
This immediately yields our next result.

\begin{proposition}
\label{prop:3}
Let $A$ be an abelian group and let $f$ be a 
$G$-equivariant 
function from the set of subgroups of $G$ to $A$. 
Then
\[
\sum_{C\in \CR/G}(-1)^{\Size{C}} f(G_C) =
\sum_{C\in \CP/G}(-1)^{\Size{C}} f(G_C).
\]
\end{proposition}

Let $\CB$ denote the spherical Tits building of $G$. 
We can view $\CB$ as the complex consisting of the parabolic subgroups
of $G$ with reversed inclusion giving the poset structure, \cite{tits}.
Then, as complexes, $\CP$ is simply the barycentric subdivision of $\CB$
and thus both $\CP$ and $\CB$ are homotopy equivalent,  
cf.\ \cite[(1.4)]{quillen}.
Let $T$ be a maximal torus of $G$ and let $S$
be a set of simple roots of $G$ with respect to $T$.
For a subset $I$ of $S$ let 
$P_I$ be the (standard) parabolic subgroup of $G$ 
associated with $I$ and 
note any parabolic subgroup of $G$ is conjugate to $P_I$ 
for some $I \subseteq S$.

Our next result follows from Proposition \ref{prop:3}, the 
comments in the previous paragraph, and the fact
that a parabolic subgroup of $G$ is self-normalizing.
For a parabolic subgroup $P$ of $G$ 
conjugate to $P_I$
let $\cor(P)$ denote its (semisimple) \emph{corank} in $G$, 
i.e.\ $\cor(P_I) = \Size{S\setminus I}$.

\begin{proposition}
\label{prop:2}
Let $A$ be an abelian group and let $f$ be a 
$G$-equivariant 
function from the set of subgroups of $G$ to $A$. 
Then
\[
\sum_{C\in \CP/G}(-1)^{\Size{C}} f(G_C) =
\sum_{P \in \CB/G}(-1)^{\cor(P)} f(P) =
\sum_{I \subseteq S}(-1)^{\Size{S\setminus I}} f(P_I).
\]
\end{proposition}

\begin{remark}
Let $p > 0$ be the characteristic of $k$ and let $q = p^a$ for some
integer $a$.
For a suitable choice of a $G$-equivariant 
function $f$ in Proposition \ref{prop:2},
the resulting sum gives the number of isomorphism classes of 
projective simple $kG(q)$-modules, where $G(q)$ is the finite group of 
Lie Type associated to $G$ and $q$, \cite[Thm.\ 5.3]{knoerrrobinson}. 
\end{remark}


\begin{remark}
There are analogues of all the results above, using 
complexes of chains of unipotent subgroups of $G$ 
(with $G$ acting by conjugation) in place of nilpotent
subalgebras of $\fg$.
We leave the details to the reader.
\end{remark}

\begin{remark}
All results above are independent of the characteristic of 
the underlying field. In particular, they are valid even if 
the characteristic of $k$ is a \emph{bad} prime for $G$, leading to 
degeneracies in the commutator relations.
Amusingly, all the results above are also true in characteristic zero.
\end{remark}

\bigskip

{\bf Acknowledgements}: 
We would like to thank G.R.\ Robinson for helpful discussions.


\end{document}